\documentclass[12pt]{amsart}
\usepackage{amsmath,amssymb,graphicx,enumerate,cite}

\usepackage[varg]{txfonts}
\newcommand\BS\boldsymbol
\newcommand\dif{\mathrm{d}}

\newcommand\expec{\mathbb E}
\newcommand\expect[4]{\expec_{#2}^{#3,#4}\left[#1\right]}
\newcommand\expecta[5]{\expec_{#2,#5}^{#3,#4}\left[#1\right]}
\newcommand\kolm[2]{L_{FP}(#1,#2)}
\newcommand\kolma[3]{L_{FP,#3}(#1,#2)}

\newcommand\Kolm[3]{\mathcal L_{K}(#1,#2,#3)}
\newcommand\Kolma[4]{\mathcal L_{K,#4}(#1,#2,#3)}
\newcommand\Kolmd[3]{\mathcal L_{K}^\dagger(#1,#2,#3)}
\newcommand\Kolmad[4]{\mathcal L_{K,#4}^\dagger(#1,#2,#3)}
\newcommand\soln[2]{\phi^{#1,#2}}
\newcommand\solna[3]{\phi_{#3}^{#1,#2}}
\newcommand\texp{\mathcal T\!\!\exp}
\newcommand\tmap[3]{\BS T_{#1}^{#2,#3}}
\newcommand\tmapp[5]{\BS T_{#1}^{#2,#3}}

\newcommand\picturehere[1]{\resizebox{0.5\textwidth}{!}{\includegraphics{#1}}}
\newcommand\parderiv[2]{\frac{\partial #1}{\partial #2}}
\newenvironment{acknowledgment}{\vspace{1EM}{\bf Acknowledgment.}}{}

\numberwithin{equation}{section}

\begin{document}

\title{Improved linear response for stochastically driven systems}

\author{Rafail V. Abramov}

\address{Department of Mathematics, Statistics and Computer
  Science\\University of Illinois at Chicago\\851 S. Morgan st. (M/C
  249)\\ Chicago, IL 60607}

\email{abramov@math.uic.edu}


\date{\today}

\begin{abstract}
The recently developed short-time linear response algorithm, which
predicts the average response of a nonlinear chaotic system with
forcing and dissipation to small external perturbation, generally
yields high precision of the response prediction, although suffers
from numerical instability for long response times due to positive
Lyapunov exponents. However, in the case of stochastically driven
dynamics, one typically resorts to the classical
fluctuation-dissipation formula, which has the drawback of explicitly
requiring the probability density of the statistical state together
with its derivative for computation, which might not be available with
sufficient precision in the case of complex dynamics (usually a
Gaussian approximation is used). Here we adapt the short-time linear
response formula for stochastically driven dynamics, and observe that,
for short and moderate response times before numerical instability
develops, it is generally superior to the classical formula with
Gaussian approximation for both the additive and multiplicative
stochastic forcing. Additionally, a suitable blending with classical
formula for longer response times eliminates numerical instability and
provides an improved response prediction even for long response times.
\end{abstract}

\maketitle

\section{Introduction}

The fluctuation-dissipation theorem (FDT) is one of the cornerstones
of modern statistical physics. Roughly speaking, the
fluctuation-dissipation theorem states that for dynamical systems at
statistical equilibrium the average response to small external
perturbations can be calculated through the knowledge of suitable
correlation functions of the unperturbed dynamical system. The
fluctuation-dissipation theorem has great practical use in a variety
of settings involving statistical equilibrium of baths of identical
gas or liquid molecules, Ornstein-Uhlenbeck Brownian motion, motion of
electric charges, turbulence, quantum field theory, chemical physics,
physical chemistry and other areas. The general advantage provided by
the fluctuation-dissipation theorem is that one can successfully
predict the response of a dynamical system at statistical equilibrium
to an arbitrary small external perturbation without ever observing the
behavior of the perturbed system, which offers great versatility and
insight in understanding behavior of dynamical processes near
equilibrium in numerous scientific applications
\cite{EvaMor,KubTodHas}. In particular, there has been a profound
interest among the atmospheric/ocean science community to apply the
fluctuation-dissipation theorem to predict global climate changes
responding to variation of certain physical parameters
\cite{Bel,Lei,CarFalIsoPurVul,GriBra,Gri,GriBraDym,GriBraMaj,GriDym,MajAbrGro,CohCra},
where the FDT has been used largely in its classical formulation
\cite{Ris}. A vivid demonstration of high predictive skill in
low-frequency climate response despite structural instability of
statistical states is given in \cite{MajAbrGer}.

Recently, Majda and the author \cite{AbrMaj5,AbrMaj4,AbrMaj6}
developed and tested a novel computational algorithm for predicting
the mean response of nonlinear functions of states of a chaotic
dynamical system to small change in external forcing based on the
FDT. The major difficulty in this situation is that the probability
measure in the limit as time approaches infinity in this case is
typically a Sinai-Ruelle-Bowen probability measure which is supported
on a large-dimensional (often fractal) set and is usually not
absolutely continuous with respect to the Lebesgue measure
\cite{EckRue,You}. In the context of Axiom A attractors, Ruelle
\cite{Rue1,Rue2} has adapted the classical calculations for FDT to
this setting. The geometric algorithm (also called the short-time FDT,
or ST-FDT algorithm in \cite{AbrMaj5,AbrMaj4,AbrMaj6}) is based on the
ideas of \cite{Rue,Rue2} and takes into account the fact that the
dynamics of chaotic nonlinear forced-dissipative systems often reside
on chaotic fractal attractors, where the classical FDT formula of the
fluctuation-dissipation theorem often fails to produce satisfactory
response prediction, especially in dynamical regimes with weak and
moderate chaos and slower mixing. It has been discovered in
\cite{AbrMaj4,AbrMaj5,AbrMaj6} that the ST-FDT algorithm is an
extremely precise response approximation for short response times, and
can be blended with the classical FDT algorithm with Gaussian
approximation of the state probability density (quasi-Gaussian FDT
algorithm, or qG-FDT) for longer response times to alleviate
undesirable effects of expanding Lyapunov directions (which cause
numerical instability in ST-FDT for longer response times). Further
developing the ST-FDT response algorithm for practical applications,
in \cite{Abr5} the author designed a computationally inexpensive
method for ST-FDT using the reduced-rank tangent map, and in
\cite{Abr6} the ST-FDT algorithm is adapted for the response on slow
variables of multiscale dynamics, which improves its computational
stability and simultaneously reduces computational expense.

However, dynamical systems describing real-world processes are often
driven by a stochastic forcing. In this setting, the traditional
approach is to use the classical FDT algorithm, which computes the
linear response to small external forcing as a correlation function
along a single long-term trajectory. Typically, it is assumed that the
single long-term trajectory samples the statistical equilibrium state
of the model, however, suitable generalizations for dynamics with
time-periodic forcing can also be made \cite{MajWan,MajGer}. A
significant drawback of the classical FDT approach is that its
computational algorithm requires the statistical state probability
density together with its derivative to be explicitly computed, which
is typically not possible for complex nonlinear systems. Usually, an
approximation is used, such as the Gaussian approximation with
suitable mean state and covariance matrix
\cite{AbrMaj4,AbrMaj5,AbrMaj6}. In this case, if the actual
statistical state is far from the Gaussian, the predicted response is
usually considerably different from what is observed by direct model
perturbation (so called ideal response
\cite{AbrMaj4,AbrMaj5,AbrMaj6}).

On the other hand, the ST-FDT response algorithm is observed to be
consistently superior to the classical FDT with Gaussian approximation
for deterministic chaotic dynamical systems with strongly non-Gaussian
statistical states for response times before the numerical instability
occurs. In this work we adapt the ST-FDT linear response algorithm to
be used with stochastically forced dynamics (further called stochastic
ST-FDT, or SST-FDT). Below we observe that the SST-FDT response
algorithm, adapted to stochastically driven dynamics and blended with
the qG-FDT algorithm to avoid numerical instability, is also generally
superior to the classical FDT with Gaussian approximation of the
statistical state for both the additive and multiplicative noise, just
as the ST-FDT algorithm in \cite{AbrMaj4,AbrMaj5,AbrMaj6} for chaotic
deterministic systems. The manuscript is organized as follows. In
Section \ref{sec:fdt_stoch} we develop the SST-FDT formula for general
time-dependent stochastically forced dynamics, and design a practical
computational algorithm for autonomous dynamics with invariant
probability measure. In Section \ref{sec:l96_app} we test the new
algorithm for the stochastically driven Lorenz 96 model
\cite{Lor,LorEma}. Section \ref{sec:summary} summarizes the results of
this work.

\section{Fluctuation-dissipation theorem for stochastically driven
systems}
\label{sec:fdt_stoch}

Here we consider an It\=o stochastic differential equation (SDE) of the
form
\begin{equation}
\label{eq:dyn_sys}
\dif\BS x = \BS f_\alpha(\BS x,t)\dif t+\BS\sigma(\BS x,t)\dif\BS W_t,
\end{equation}
where $\BS x=\BS x(t)\in\mathbb R^N$, $\BS f_\alpha:[\mathbb R^N\times
  T]\to\mathbb R^N$, $\BS \sigma:[\mathbb R^{N\times K}\times
  T]\to\mathbb R^N$ are smooth nonlinear functions, and $\BS W_t$ is the
$K$-dimensional Wiener process. Additionally, $\BS f$ depends on a
scalar parameter $\alpha$. We say that the SDE in \eqref{eq:dyn_sys}
is {\em unperturbed} if $\alpha=0$, or {\em perturbed} otherwise. We
also adopt the notation $\BS f\equiv\BS f_0$, with the
assumption
\begin{equation}
\label{eq:f_alpha}
\parderiv{}\alpha\BS f_\alpha(\BS x,t)|_{\alpha=0}=\BS B(\BS x)\BS\eta(t),
\end{equation}
where $\BS B(\BS x)$ is an $N\times L$ matrix-valued function, and
$\BS \eta(t)$ is a $L$-vector valued function. The practical meaning
of the above assumption will become clear below.

Let $A(\BS x)$ be a nonlinear function of $\BS x$, and let
$\expecta{A}{\BS x}{t_0}{t}{\alpha}$, where $t>0$ is the elapsed time
after $t_0$, denote the expectation of $A$ at time $t_0+t$ over all
realizations of the Wiener process in \eqref{eq:dyn_sys}, under the
condition that $\BS x(t_0)=\BS x$ (with the short notation
$\expecta{A}{\BS x}{t_0}{t}{0}=\expect{A}{\BS x}{t_0}{t}$). Let $A$ at
the time $t_0$ be distributed according to a probability measure
$\rho_{t_0}$, that is, the average value of $A$ at time $t_0$ is
\begin{equation}
\langle A\rangle(t_0)=\rho_{t_0}(A)=\int_{\mathbb R^N} A(\BS
x)\dif\rho_{t_0}(\BS x),
\end{equation}
where $\dif\rho_{t_0}(\BS x)$ denotes the measure of the infinitesimal
Lebesgue volume $\dif\BS x$ associated with $\BS x$. Then, for time
$t_0+t$, the average of $A$ for the perturbed system in
\eqref{eq:dyn_sys} is given by
\begin{equation}
\langle A\rangle_\alpha(t_0+t)=\rho_{t_0}\left(\expecta{A}{\BS x}
{t_0}{t}{\alpha}\right)=\int_{\mathbb R^N}
\expecta{A}{\BS x}{t_0}{t}{\alpha}\dif\rho_{t_0}(\BS x).
\end{equation}
In general, for the same initial distribution $\rho_{t_0}$, the
average value $\langle A\rangle_\alpha(t_0+t)$ depends on the value of
$\alpha$. Here, we define the {\em average response} $\delta\langle
A\rangle_\alpha(t_0+t)$ as
\begin{equation}
\label{eq:av_resp}
\delta\langle A\rangle_\alpha(t_0+t)=\int_{\mathbb
R^N}\left(\expecta{A}{\BS x}{t_0}{t}{\alpha}-
\expect{A}{\BS x}{t_0}{t}\right)\dif\rho_{t_0}(\BS x).
\end{equation}
The meaning of the average response in \eqref{eq:av_resp} is the
following: for the same initial average value of $A$ it provides the
difference between the future average values of $A$ for the perturbed
and unperturbed dynamics in \eqref{eq:dyn_sys}.

If $\alpha$ is small, we can formally linearize \eqref{eq:av_resp}
with respect to $\alpha$ by expanding in Taylor series around
$\alpha=0$ and truncating to the first order, obtaining the following
general linear fluctuation-response formula:
\begin{equation}
\label{eq:lin_resp}
\delta\langle A\rangle_\alpha(t_0+t)=\alpha\int_{\mathbb R^N}\partial_\alpha
\expect{A}{\BS x}{t_0}{t}\dif\rho_{t_0}(\BS x),
\end{equation}
where we use the short notation
\begin{equation}
\partial_\alpha\bullet\equiv\left.\parderiv{\bullet_\alpha}
\alpha\right|_{\alpha=0}.
\end{equation}
\subsection{Stochastic short-time linear response}
To compute the general linear fluctuation-response formula in
\eqref{eq:lin_resp}, we need a suitable algorithm for
$\partial_\alpha\expect{A}{\BS x}{t_0}{t}$. Let
$x(t_0+t)=\solna{t_0}{t}{\alpha}\BS x$ be the trajectory of
\eqref{eq:dyn_sys} starting at $\BS x$ at $t_0$ for a particular
realization of the Wiener process $\BS W_{[t_0\ldots t_0+t]}$. Then, the
expectation $\expecta{A}{\BS x}{t_0}{t}{\alpha}$ is given by
\begin{equation}
\expecta{A}{\BS x}{t_0}{t}{\alpha}=\expec\left[A
\left(\solna{t_0}{t}{\alpha}\BS x\right)\right],
\end{equation}
where the expectation in the right-hand side is taken with respect to
all Wiener paths. Therefore,
\begin{equation}
\label{eq:expect_partial}
\partial_\alpha\expect{A}{\BS x}{t_0}{t}=
\expec\left[DA\left(\soln{t_0}{t}\BS x\right)
\partial_\alpha\soln{t_0}{t}\BS x\right],
\end{equation}
where $DA$ denotes the derivative of $A$ with respect to its
argument. For $\partial_\alpha\soln{t_0}{t}\BS x$, by taking the
difference between the perturbed and unperturbed versions of
\eqref{eq:dyn_sys} and linearizing with respect to $\alpha$ at
$\alpha=0$, we have
\begin{equation}
\label{eq:dif_soln}
\begin{split}
\dif\partial_\alpha\soln{t_0}{t}\BS x=\Big(D\BS
f\left(\soln{t_0}{t}\BS x, t_0+t\right)\dif t+D\BS
\sigma\left(\soln{t_0}{t}\BS x,t_0+t\right) \dif\BS
W_{t_0+t}\Big)\times\\\times
\partial_\alpha\soln{t_0}{t}\BS x+\partial_\alpha\BS f
\left(\soln{t_0}{t}\BS x,t_0+t\right)\dif t,
\end{split}
\end{equation}
where $D\BS f$ and $D\BS \sigma$ are Jacobians of $\BS f$ and
$\BS\sigma$, respectively. The above equation is a linear stochastic
differential equation for $\partial_\alpha\soln{t_0}{t}\BS x$ with
zero initial condition (as at $t_0$ both perturbed and unperturbed
solutions start with the same $\BS x$). It can be solved as follows:
let us first introduce the integrating factor $\tmap{\BS x}{t_0}{t}$
(an $N\times N$ matrix) given by the solution of the equation
\begin{equation}
\label{eq:tmap}
\begin{split}
\dif\tmap{\BS x}{t_0}{t}&=\Big(D\BS f\left(\soln{t_0}{t}\BS
x,t_0+t\right) \dif t+\\&+D\BS\sigma\left(\soln{t_0}{t}\BS
x,t_0+t\right)\dif\BS W_{t_0+t}\Big) \tmap{\BS x}{t_0}{t},\quad
\tmap{\BS x}{t_0}{0}=\BS I,
\end{split}
\end{equation}
and represent $\partial_\alpha\soln{t_0}{t}\BS x$ as a product
\begin{equation}
\partial_\alpha\soln{t_0}{t}\BS x=\tmap{\BS x}{t_0}{t}\BS
y^{t_0,t}_{\BS x},
\end{equation}
where $\BS y^{t_0,t}_x$ is an $N$-vector. Then, for the It\=o
differential of $\partial_\alpha\soln{t_0}{t}\BS x$ we obtain
\begin{equation}
\label{eq:dif_soln_2}
\begin{split}
\dif\partial_\alpha\soln{t_0}{t}\BS x=\dif\tmap{\BS x}{t_0}{t}\BS
y^{t_0,t}_{\BS x}+ \tmap{\BS x}{t_0}{t}\dif\BS y^{t_0,t}_{\BS
  x}=\Big(D\BS f\left(\soln{t_0}{t}\BS x,t_0+t\right) \dif
t+\\+D\BS\sigma\left(\soln{t_0}{t}\BS x,t_0+t\right)\dif\BS
W_{t_0+t}\Big) \tmap{\BS x}{t_0}{t}\BS y^{t_0,t}_{\BS x}+\tmap{\BS
  x}{t_0}{t}\dif\BS y^{t_0,t}_{\BS x}=\\ =\Big(D\BS
f\left(\soln{t_0}{t}\BS x,t_0+t\right)\dif
t+D\BS\sigma\left(\soln{t_0}{t}\BS x,t_0+t\right)\dif\BS
W_{t_0+t}\Big)\times\\ \times\partial_\alpha\soln{t_0}{t}\BS
x+\tmap{\BS x}{t_0}{t}\dif\BS y^{t_0,t}_{\BS x}.
\end{split}
\end{equation}
Comparing the right-hand sides of \eqref{eq:dif_soln} and
\eqref{eq:dif_soln_2} we find that $\BS y^{t_0,t}_{\BS x}$ satisfies
\begin{equation}
\dif\BS y^{t_0,t}_{\BS x}=(\tmap{\BS
x}{t_0}{t})^{-1}\partial_\alpha\BS f \left(\soln{t_0}{t}\BS
x,t_0+t\right)\dif t,\quad\BS y^{t_0,0}_{\BS x}=0,
\end{equation}
with the formal solution
\begin{equation}
\BS y^{t_0,t}_{\BS x}=\int_0^t(\tmap{\BS
  x}{t_0}{\tau})^{-1}\partial_\alpha\BS f
\left(\soln{t_0}{\tau}\BS x,t_0+\tau\right)\dif\tau.
\end{equation}
Therefore, $\partial_\alpha\soln{t_0}{t}\BS x$ is given by
\begin{equation}
\label{eq:dif_soln_3}
\partial_\alpha\soln{t_0}{t}\BS x= \int_0^t\tmap{\BS
  x}{t_0}{t}(\tmap{\BS x}{t_0}{\tau})^{-1}\partial_\alpha\BS f
\left(\soln{t_0}{\tau}\BS x,t_0+\tau\right)\dif\tau.
\end{equation}
At this point, observe that the solution $\tmap{\BS x}{t_0}{t}$ of
\eqref{eq:tmap} can be represented as a product
\begin{equation}
\tmap{\BS x}{t_0}{t}=\tmap{\soln{t_0}{\tau}\BS x}{t_0+\tau}
{t-\tau}\tmap{\BS x}{t_0}{\tau},\quad \tau\leq t,
\end{equation}
due to the fact that a solution of \eqref{eq:tmap} can be multiplied
by an arbitrary constant matrix on the right and still remains the
solution. Then, \eqref{eq:dif_soln_3} becomes
\begin{equation}
\label{eq:deriv_soln}
\begin{split}
\partial_\alpha\soln{t_0}{t}\BS x=\int_0^t\tmapp{\soln{t_0}{\tau}\BS x}
{t_0+\tau}{t-\tau}{t_0+\tau}{t_0+t}\partial_\alpha\BS f
\left(\soln{t_0}{\tau}\BS x,t_0+\tau\right)\dif\tau.
\end{split}
\end{equation} 
For smooth $\BS f_\alpha$ and $\BS \sigma$ in \eqref{eq:dyn_sys},
$\soln{t_0}{t}\BS x$ smoothly depends on $\BS x$ \cite{Kun}, and the
integrating factor $\tmap{\BS x}{t_0}{t}$ is in fact the tangent map for
the trajectory $\soln{t_0}{t}\BS x$:
\begin{equation}
\tmap{\BS x}{t_0}{t}=\parderiv{}{\BS x}\soln{t_0}{t}\BS x.
\end{equation}
With \eqref{eq:deriv_soln},
\eqref{eq:expect_partial} becomes
\begin{equation}
\partial_\alpha\expect{A}{\BS x}{t_0}{t}=\int_0^t\expec
\Big[DA\left(\soln{t_0}{t}\BS x\right)\tmapp{\soln{t_0}{\tau}\BS x}
{t_0+\tau}{t-\tau}{t_0+\tau}{t_0+t}
\partial_\alpha\BS f\left(\soln{t_0}{\tau}\BS x,t_0+\tau\right)
\Big]\dif\tau.
\end{equation}
Recalling \eqref{eq:f_alpha}, we write the above formula as
\begin{equation}
\partial_\alpha\expect{A}{\BS x}{t_0}{t}=\int_0^t
\expec\Big[DA\left(\soln{t_0}{t}\BS x\right)
\tmapp{\soln{t_0}{\tau}\BS x}{t_0+\tau}{t-\tau}{t_0+\tau}
{t_0+t}\BS B\left(\soln{t_0}{\tau}\BS x\right)\Big]\BS
\eta(t_0+\tau)\dif\tau.
\end{equation}
Then, the general linear response formula in \eqref{eq:lin_resp} can be
written as
\begin{equation}
\label{eq:fdt_response}
\delta\langle A\rangle_\alpha(t_0+t)=\alpha\int_0^t
\BS R_{SST}(t_0,t,\tau)\BS \eta(t_0+\tau)\dif\tau,
\end{equation}
where the {\em linear response operator} $\BS R_{SST}(t_0,t,\tau)$ is given by
\begin{equation}
\label{eq:fdt_operator}
\BS R_{SST}(t_0,t,\tau)=\expec\int_{\mathbb
R^N}DA\left(\soln{t_0}{t}\BS x\right)
\tmapp{\soln{t_0}{\tau}\BS x}{t_0+\tau}{t-\tau}{t_0+\tau}{t_0+t}
\BS B\left(\soln{t_0}{\tau}\BS x\right)\dif\rho_{t_0}(\BS x).
\end{equation}
Further we refer to \eqref{eq:fdt_operator} as the {\em stochastic
  short-time fluctuation-dissipation theorem algorithm}, or SST-FDT
algorithm. The reason is that in practice the computation of the
tangent map in \eqref{eq:tmap} for large $t$ becomes numerically
unstable because of exponential growth due to positive Lyapunov
exponents (just as observed in
\cite{Abr5,Abr6,AbrMaj4,AbrMaj5,AbrMaj6} for deterministic chaotic
dynamics). Note that if the stochastic forcing is removed from
\eqref{eq:dyn_sys}, the SST-FDT response operator becomes the usual
ST-FDT from \cite{Abr5,Abr6,AbrMaj4,AbrMaj5,AbrMaj6}. Apparently,
\eqref{eq:fdt_operator} requires the average with respect to
$\rho_{t_0}$. If $\rho_{t_0}$ is not known explicitly, there are some
opportunities to replace the $\rho$-average with time average,
particularly for the autonomous dynamics with $\rho_{t_0}$ being the
invariant probability measure, and also for non-autonomous dynamical
systems with explicit time-periodic dependence (as done in
\cite{MajWan,MajGer} for classical FDT response).

\subsection{Classical linear response}

The standard way to derive the classical linear response formula is
through the Fokker-Planck equation (or, as it is also called, the
forward Kolmogorov equation) for the perturbed system in
\eqref{eq:dyn_sys} by neglecting the terms of higher order than the
perturbation, as it is done in
\cite{AbrMaj4,AbrMaj5,AbrMaj6,MajAbrGro,Ris,MajWan}. However, for the
sake of clarity, here we show the derivation of the classical FDT
directly from \eqref{eq:lin_resp}. Under the assumption of continuity
of $\rho_{t_0}$ with respect to the Lebesgue measure, that is,
$\dif\rho_{t_0}(\BS x)=p_{t_0}(\BS x)\dif\BS x$, where $p_{t_0}$ is the
probability density, we can also obtain a formal general expression for
the classical fluctuation-response formula. Using the notations
\begin{equation}
\label{eq:kolm-notations}
\begin{array}{c}
\displaystyle\kolma{\BS x}{t}{\alpha}=-\parderiv{}{\BS x}\cdot
(\BS f_\alpha(\BS x,t)\bullet)+\left(\parderiv{}{\BS x}
\otimes\parderiv{}{\BS x}\right)
\cdot(\BS{\sigma\sigma}^T(\BS x,t)\bullet),\\
\displaystyle\Kolma{\BS x}{t_0}{t}{\alpha}=
\texp\left(\int_0^t\dif\tau \,\kolma{\BS x}{t_0+\tau}{\alpha}\right),
\end{array}
\end{equation}
which are, respectively, the Fokker-Planck and forward Kolmogorov
operators, we write the expectation $\expecta{A}{\BS x}{t_0}{t}{\alpha}$
in the form
\begin{equation}
\begin{split}
&\expecta{A}{\BS x}{t_0}{t}{\alpha}=\int_{\mathbb R^N}A(\BS
y)\Kolma{\BS y}{t_0}{t}{\alpha}\delta(\BS x-\BS y)\dif\BS
y=\\ =\int_{\mathbb R^N}&\Kolma{\BS y}{t_0}{t}{\alpha}A(\BS
y)\delta(\BS x-\BS y) \dif\BS y=\Kolmad{\BS x}{t_0}{t}{\alpha}A(\BS
x),
\end{split}
\end{equation}
where $\delta(\BS x)$ is the Dirac delta-function, and the adjoint is
taken with respect to the standard inner product under the
integral. Then, the general response formula with $\dif\rho_{t_0}(\BS
x)= p_{t_0}(\BS x)\dif\BS x$ becomes
\begin{equation}
\begin{split}
\delta\langle A\rangle_\alpha(t_0+t)&=\alpha\int_{\mathbb
  R^N}\partial_\alpha \expect{A}{\BS x}{t_0}{t}p_{t_0}(\BS x)\dif \BS
x=\\=&\alpha\int_{\mathbb R^N}\partial_\alpha \Kolmd{\BS
  x}{t_0}{t}A(\BS x)p_{t_0}(\BS x)\dif \BS x=\\=&\alpha\int_{\mathbb
  R^N}A(\BS x)\partial_\alpha\Kolm{\BS x}{t_0}{t}p_{t_0}(\BS x)\dif\BS
x.
\end{split}
\end{equation}
It is not difficult to show that the parametric derivative of an
ordered exponential of a linear operator $L_\alpha(\BS x,t)$ is
computed as
\begin{equation}
\label{eq:texp_deriv}
\begin{split}
\parderiv{}\alpha&\texp\left(\int_{t_0}^t\dif\tau\,L_\alpha(\BS
x,\tau)\right)=\\ =&\int_{t_0}^t\dif\tau\,\texp\left(\int_\tau^t\dif
s\,L_\alpha(\BS x,s)\right) \parderiv{L_\alpha(\BS x,\tau)}\alpha\texp
\left(\int_{t_0}^\tau\dif s\,L_\alpha(\BS x,s)\right).
\end{split}
\end{equation}
As a result, we obtain
\begin{equation}
\begin{split}
\delta\langle A&\rangle_\alpha(t_0+t)=
\alpha\int_0^t\dif\tau\,\int_{\mathbb R^N}\Kolmd{\BS
  x}{t_0+\tau}{t-\tau} \times\\\times &A(\BS
x)\partial_\alpha\kolm{\BS x}{t_0+\tau}\Kolm{\BS x}{t_0}{\tau}
p_{t_0}(\BS x)\dif\BS x=\\=\alpha\int_0^t&\dif\tau\,\int_{\mathbb R^N}
\expect{A}{\BS x}{t_0+\tau}{t-\tau}\partial_\alpha\kolm{\BS
  x}{t_0+\tau} p_{t_0+\tau}(\BS x)\dif\BS x,
\end{split}
\end{equation}
where $p_{t_0+\tau}(\BS x)$ is given by
\begin{equation}
p_{t_0+\tau}(\BS x)=\Kolm{\BS x}{t_0}{\tau}p_{t_0}(\BS x).
\end{equation}
Recalling \eqref{eq:f_alpha}, we recover the classical linear
fluctuation-response formula in the form
\begin{equation}
\label{eq:class_resp}
\begin{split}
\delta\langle A&\rangle_\alpha(t_0+t)=\alpha\int_0^t\dif\tau\,
\int_{\mathbb R^N}\expect{A}{\BS
  x}{t_0+\tau}{t-\tau}\times\\\times&\partial_\alpha
\kolm{\BS x}{t_0+\tau}p_{t_0+\tau}(\BS x)\dif\BS x=
\alpha\int_0^t\BS R_{class}(t_0,t,\tau)\BS\eta(t_0+\tau)\dif\tau,
\end{split}
\end{equation}
where the classical linear response operator $\BS R_{class}$ is given
by
\begin{equation}
\label{eq:class_operator}
\BS R_{class}(t_0,t,\tau)=-\expec\int_{\mathbb
  R^N}A(\soln{t_0+\tau}{t-\tau}\BS x)
\parderiv{}{\BS x}\cdot(\BS B(\BS x)p_{t_0+\tau}(\BS x))\dif\BS x.
\end{equation}
Observe that, unlike \eqref{eq:fdt_operator}, in
\eqref{eq:class_operator} one has to know $p_{t_0+\tau}(\BS x)$ for
all response times explicitly to perform differentiation with respect
to $\BS x$. Usually, an approximation is used, such as the Gaussian
approximation \cite{AbrMaj4,AbrMaj5,AbrMaj6}.

\subsection{Special case for autonomous dynamics with ergodic
invariant probability measure}

Here we consider the case where $\BS f$ and $\BS \sigma$ in
\eqref{eq:dyn_sys} do not explicitly depend on $t$ (although $\BS
f_\alpha$ does with $\alpha\neq 0$), and we choose $\rho_{t_0}=\rho$
to be an ergodic invariant probability measure for
\eqref{eq:dyn_sys}. In this situation, one can replace the averaging
with respect to the measure $\rho$ with averaging over a single
long-term trajectory which starts with an initial condition $\BS x$ in
the support of $\rho$:
\begin{equation}
\begin{split}
\BS R_{SST}(t_0,t,\tau)=\expec\lim_{r\to\infty}\frac 1r\int_0^rDA
\left(\soln{t_0}{t}\soln{t_0-s}{s}\BS x\right)\times\\\times
\tmapp{\soln{t_0}{\tau}\soln{t_0-s}{s}\BS
x}{t_0+\tau}{t-\tau}{t_0+\tau}{t_0+t}
\BS B\left(\soln{t_0}{\tau}\soln{t_0-s}{s}\BS x\right)\dif s,
\end{split}
\end{equation}
where, without loss of generality, the starting is time $t_0-s$, that
is, the averaging occurs over the endpoints of $\soln{t_0-s}{s}\BS
x$. Combining the solution operators, we obtain
\begin{equation}
\begin{split}
\BS R_{SST}(t_0,t,\tau)=\expec\lim_{r\to\infty}\frac 1r\int_0^rDA
\left(\soln{t_0-s}{s+t}\BS x\right)\times\\\times
\tmapp{\soln{t_0-s}{s+\tau}\BS x}{s+\tau}{t-\tau}{t_0+\tau}{t_0+t}
\BS B\left(\soln{t_0-s}{s+\tau}\BS x\right)\dif s.
\end{split}
\end{equation}
Since the averaging over all independent realizations of the Wiener
process is needed, we can average over many statistically independent
chunks of the Wiener path along a single long-time trajectory by
setting $t_0=s-\tau$:
\begin{equation}
\BS R_{SST}(t,\tau)=\lim_{r\to\infty}\frac 1r\int_0^rDA
\left(\soln{-\tau}{s+t}\BS x\right)
\tmapp{\soln{-\tau}{s+\tau}\BS x}{s}{t-\tau}{t_0+\tau}{t_0+t}
\BS B\left(\soln{-\tau}{s+\tau}\BS x\right)\dif s.
\end{equation}
Finally, replacing $\BS x$ with $\soln{0}{-\tau}\BS x$ (which for
finite $\tau$ is also in the support of $\rho$), we find that
\begin{equation}
\BS R_{SST}(t,\tau)=\lim_{r\to\infty}\frac 1r\int_0^rDA
\left(\soln{0}{s+t-\tau}\BS x\right)
\tmapp{\soln{0}{s}\BS x}{s}{t-\tau}{t_0+\tau}{t_0+t}
\BS B\left(\soln{0}{s}\BS x\right)\dif s,
\end{equation}
or, denoting $\BS x(s)=\soln{0}{s}\BS x$,
\begin{equation}
\BS R_{SST}(t,\tau)=\lim_{r\to\infty}\frac 1r\int_0^rDA(
\BS x(s+t-\tau))\tmapp{\BS x(s)}
{s}{t-\tau}{s}{s+t-\tau}\BS B\left(\BS x(s)\right)\dif s.
\end{equation}
Now, the linear response formula in \eqref{eq:fdt_response} and the
response operator in \eqref{eq:fdt_operator} become, respectively,
\begin{equation}
\label{eq:st_fdt}
\begin{split}
\delta\langle A\rangle&_\alpha(t_0+t)=\alpha\int_0^t \BS
R_{SST}(t-\tau)\BS\eta(t_0+\tau)\dif\tau,\\ \BS
R_{SST}(t)&=\lim_{r\to\infty}\frac 1r\int_0^rDA( \BS x(s+t))\tmapp{\BS
  x(s)}{s}{t}{s}{s+t}\BS B\left(\BS x(s)\right)\dif s.
\end{split}
\end{equation}
In a similar fashion, for the classical linear response in
\eqref{eq:class_resp} we note that the Fokker-Planck operator $L_{FP}$
does not depend on $t$, and both the forward Kolmogorov operator
$\mathcal L_K$ and its adjoint do not depend on $t_0$. Taking into
account that $p_{t_0+\tau}(\BS x)=p(\BS x)$, where $p(\BS x)$ is the
invariant probability density, we write
\begin{equation}
\BS R_{class}(t)=-\expec\int_{\mathbb R^N}A(\soln{t_0}{t}\BS x)
\parderiv{}{\BS x}\cdot(\BS B(\BS x)p(\BS x))\dif\BS x,
\end{equation}
or, after replacing the $p$-average with the average over the
long-term trajectory,
\begin{equation}
\BS R_{class}(t)=-\expec\lim_{r\to\infty}\frac 1r\int_0^r A(\BS
x(s+t))\frac{\parderiv{}{\BS x}\cdot(\BS B(\BS x(s))p(\BS
x(s)))}{p(x(s))}\dif s.
\end{equation}
Here the expectation can be removed since the averaging over different
Wiener paths will automatically occur as the long time average is
computed. As a result, we obtain
\begin{equation}
\BS R_{class}(t)=-\lim_{r\to\infty}\frac 1r\int_0^rA(\BS x(s+t))
\frac{\parderiv{}{\BS x}\cdot(\BS B(\BS x(s))p(\BS x(s)))}{p(\BS
x(s))}\dif s.
\end{equation}

\section{Application for the stochastically driven Lorenz 96 model}
\label{sec:l96_app}

The 40-mode deterministic Lorenz 96 model (L96) has been introduced by
Lorenz and Emanuel \cite{Lor,LorEma} as a simple model with large
scale features of complex nonlinear geophysical systems. The
deterministic Lorenz 96 (L96) model is given by
\begin{equation}
\label{eq:L96}
\dot X_n=X_{n-1}(X_{n+1}-X_{n-2})-X_n+F,\quad 1\leq k\leq N,
\end{equation}
with periodic boundary conditions given by $X_{n\pm N}=X_n$, where
$N=40$, and $F$ being a constant forcing parameter. The model in
\eqref{eq:L96} is designed to mimic midlatitude weather and climate
behavior (in particular Rossby waves), so periodic boundary conditions
are appropriate. It is demonstrated in Chapter 2 of \cite{MajAbrGro}
that the dynamical regime of the L96 model varies with changing the
value of constant forcing $F$: weakly chaotic dynamical regimes with
$F=5,6$, strongly chaotic regime with $F=8$, and turbulent regimes
$F=12,16,24$ with self-similar time autocorrelation decay.

Here we apply the stochastic forcing to the L96 model as
\begin{equation}
\label{eq:SL96}
\dif X_k=\left[X_{k-1}(X_{k+1}-X_{k-2})-X_k+F\right]\dif t
+(\BS\sigma(\BS X))_k(\dif\BS W_t)_k,
\end{equation}
where $\BS\sigma:\mathbb R^N\to\mathbb R^N$ is a vector-valued
function of $\BS X$, $\BS W$ is a $N$-dimensional Wiener process, and
$(\dif\BS W_t)_k$ is the $k$-th component of $\dif\BS W$ (that is,
effectively $\BS\sigma$ is a diagonal matrix multiplying the vector
$\dif\BS W$). As the stochastic Lorenz 96 (SL96) model above does not
depend explicitly on time (except for the Wiener noise), we can assume
that it has an invariant probability measure $\rho$.

In this work, we perturb the SL96 model in \eqref{eq:SL96} by a small
parameter $\alpha$ as
\begin{equation}
\label{eq:SL96_pert}
\dif X_k=\left[X_{k-1}(X_{k+1}-X_{k-2})-X_k+F+\alpha\eta_k
\right]\dif t+(\BS\sigma(\BS X))_k(\dif\BS W_t)_k,
\end{equation}
where $\BS\eta\in\mathbb R^N$ is a constant forcing vector
perturbation, which is ``turned on'' at time $t_0=0$. With the
invariant probability state $\rho$, and the perturbation given in
\eqref{eq:SL96_pert}, the general response formula in
\eqref{eq:lin_resp} becomes
\begin{equation}
\label{eq:resp_const_forc}
\begin{split}
\delta\langle A\rangle_\alpha(t)=\alpha\mathcal R(t)\BS\eta,\\
\mathcal R(t)=\int_0^t\BS R(\tau)\dif\tau,\\
\end{split}
\end{equation}
where subscripts for $\BS R$ and $\mathcal R$ are omitted as both the
SST-FDT and classical response operators apply. We also set the
observable $A(\BS x)=\BS x$, that is, the response of the mean state
is computed. As an approximation for the invariant probability density
for the classical response, we choose the Gaussian distribution with
the same mean and covariance as the actual invariant probability
measure, which are determined by averaging along the long-term time
series of unperturbed \eqref{eq:SL96}, and, thus, further call it
quasi-Gaussian FDT (qG-FDT) as in \cite{AbrMaj4,AbrMaj5,AbrMaj6}. In
this setting, the short-time and quasi-Gaussian linear response
operators become
\begin{equation}
\begin{split}
&\BS R_{SST}(t)=\lim_{r\to\infty}\frac 1r\int_0^r\tmapp{\BS
x(s)}{s}{t}{s}{s+t} \dif s,\\\BS R_{qG}(t)=&\lim_{r\to\infty}\frac
1r\int_0^r\BS x(s+t)\BS C^{-1} (\BS x(s)-\bar{\BS x})\dif s,
\end{split}
\end{equation}
where $\bar{\BS x}$ and $\BS C$ are the mean state and covariance
matrix of the long-time series of unperturbed \eqref{eq:SL96}.

\subsection{Blended SST/qG-FDT response}

Following \cite{AbrMaj5,AbrMaj6}, we also compute the blended
SST/qG-FDT response as
\begin{equation}
\label{eq:blend_fdt}
\BS R_{SST/qG}(t)=\left[1-H\left(t-t_{\mbox{\scriptsize
cutoff}}\right)\right]\BS R_{SST}(t)+H\left(t-t_{\mbox{\scriptsize
cutoff}}\right)\BS R_{qG}(t),
\end{equation}
where the blending function $H$ is the Heaviside step-function.  The
cut-off time $t_{\mbox{\scriptsize cutoff}}$ is chosen as
\begin{equation}
t_{\mbox{\scriptsize cutoff}}=\frac 3{\lambda_1},
\end{equation}
where $\lambda_1$ is the largest Lyapunov exponent (for details see
\cite{AbrMaj5,AbrMaj6}). This cut-off time allows to switch to the
$\BS R_{qG}$ just before the numerical instability occurs in $\BS
R_{SST}$, and, thus avoid the numerical instability. For constant
external forcing and the Heaviside blending step-function the blended
response operators become
\begin{equation}
\label{eq:STqG_mix_op}
\mathcal R_{SST/qG}(t)=\int_0^{t_{\mbox{\scriptsize cutoff}}}
\BS R_{SST}(\tau)\dif\tau + \int_{t_{\mbox{\scriptsize cutoff}}}^t
\BS R_{qG}(\tau)\dif\tau.
\end{equation}

\subsection{Computational experiments}

Below we perform computational experiments in the following setting:
\begin{itemize}
\item The number of variables (model size) $N=40$
\item Constant forcing $F=6$. The L96 model is observed to be weakly
  chaotic in this regime \cite{AbrMaj4,AbrMaj5,AbrMaj6,MajAbrGro}, and
  we would like to compare the responses for weakly chaotic
  deterministic dynamics and the stochastically driven dynamics
\item The tangent map $\tmap{x}{t_0}{t}$ in \eqref{eq:tmap} is
  computed in the same fashion as in
  \cite{Abr5,Abr6,AbrMaj4,AbrMaj5,AbrMaj6}
\item Forward Euler numerical scheme with time step $\Delta t=0.001$
  for both \eqref{eq:tmap} and \eqref{eq:SL96}
\item The linear response is tested for the following settings of the
  stochastic term $\sigma$:
  \begin{itemize}
  \item $\sigma_k=0$ (fully deterministic regime without stochastic
    forcing)
  \item $\sigma_k=1$ (additive noise)
  \item $\sigma_k=0.2X_k$, $\sigma_k=0.5X_k$ (multiplicative noise)
  \end{itemize}
\item We compute the linear response operators $\mathcal R_{SST}$,
  $\mathcal R_{qG}$ and $\mathcal R_{SST/qG}$, which are given by
  \eqref{eq:resp_const_forc} and \eqref{eq:STqG_mix_op}, and compare
  them with the ideal response operator $\mathcal R_{ideal}$, which is
  computed through the direct model perturbations
  \cite{Abr5,Abr6,AbrMaj4,AbrMaj5,AbrMaj6}
\item The time-averaging is done along a time series of 10000 time
  units
\item The ideal response operator $\mathcal R_{ideal}$ is computed via
  direct perturbations a 10000-member statistical ensemble
\item The comparison of the FDT response operators with the ideal
  response operator is carried out by evaluating the $L_2$ relative
  error
  \begin{equation}
    L_2\mbox{-error}=\frac{\|\mathcal R_{FDT}-\mathcal R_{ideal}\|}
    {\|\mathcal R_{ideal}\|},
  \end{equation}
  and the correlation function
  \begin{equation}
    \mbox{Corr}=\frac{(\mathcal R_{FDT},\mathcal R_{ideal})}
         {\|\mathcal R_{FDT}\|\|\mathcal R_{ideal}\|},
  \end{equation}
  where $(\cdot ,\cdot)$ denotes the standard Euclidean inner
  product. Observe that the $L_2$ error shows the general difference
  between the FDT and ideal responses, while the correlation function
  shows the extent to which the responses are collinear (that is, how
  well the location of the response is determined, without considering
  its magnitude)
\end{itemize}
\begin{figure}%
\picturehere{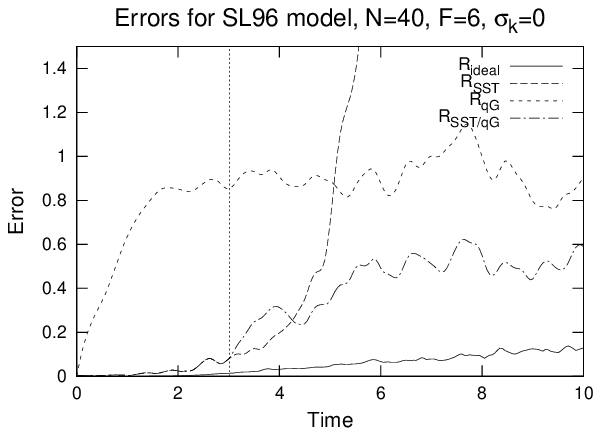}%
\picturehere{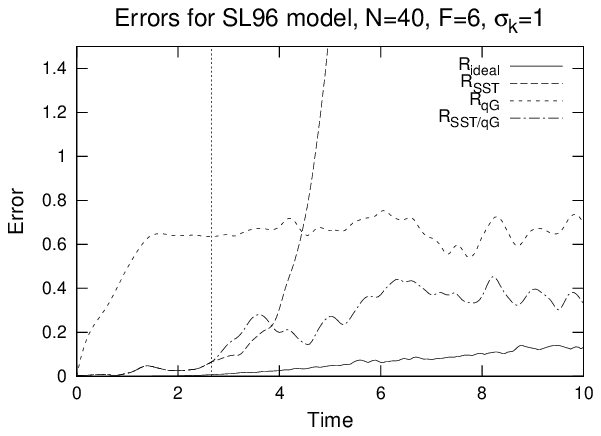}\\%
\picturehere{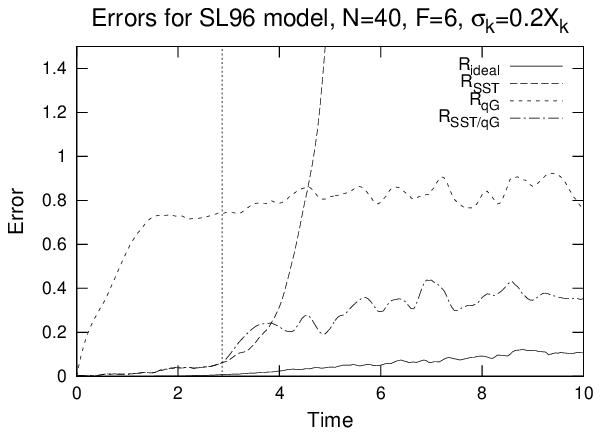}%
\picturehere{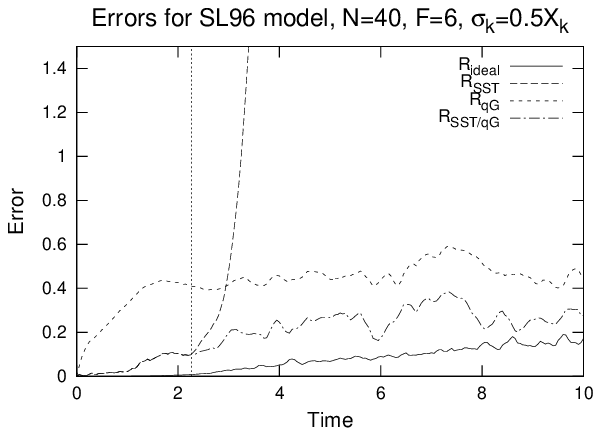}%
\caption{$L_2$-errors of the response operators for SL96 model, $N=40$,
$F=6$. Straight dotted vertical line denotes the blending cut-off time
for SST/qG-FDT. $\mathcal R_{ideal}$ denotes the intrinsic error in the ideal response due to slight nonlinearity.}%
\label{fig:l96_errors}%
\end{figure}%
\begin{figure}%
\picturehere{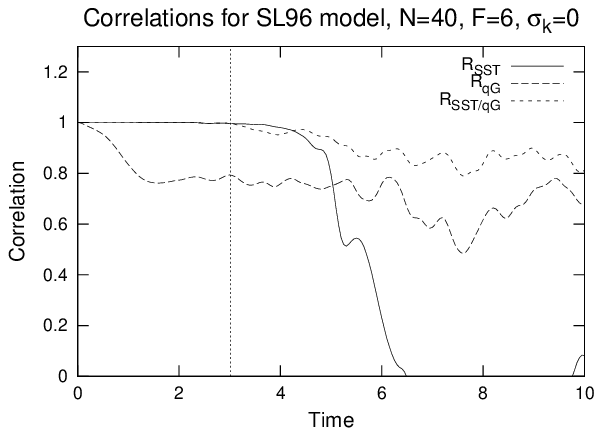}%
\picturehere{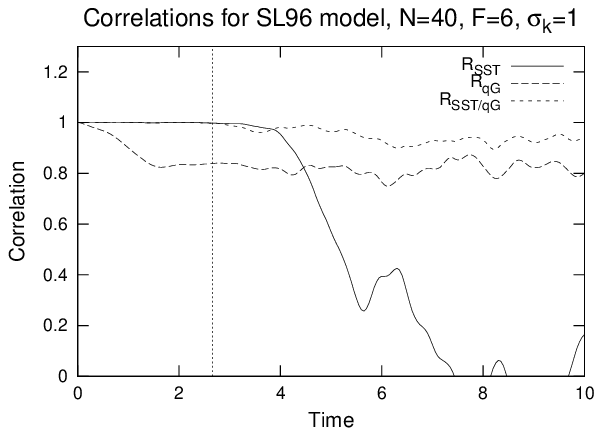}\\%
\picturehere{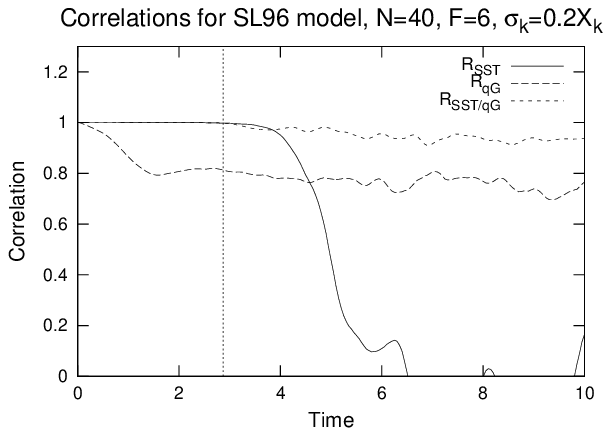}%
\picturehere{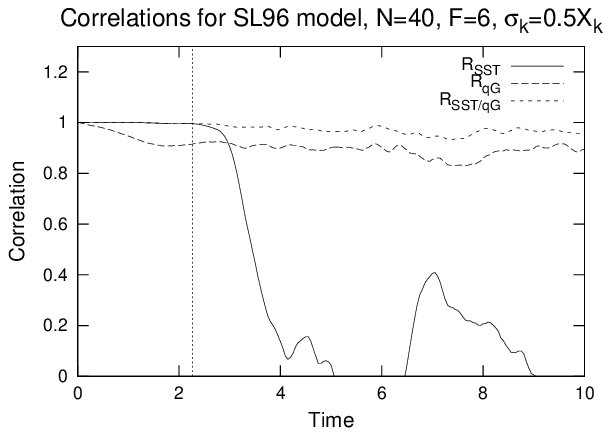}%
\caption{Correlations of the FDT response operators with the ideal
  response operator for SL96 model, $N=40$, $F=6$. Straight dotted
  vertical line denotes the blending cut-off time for SST/qG-FDT.}%
\label{fig:l96_corrs}%
\end{figure}%
\begin{figure}%
\picturehere{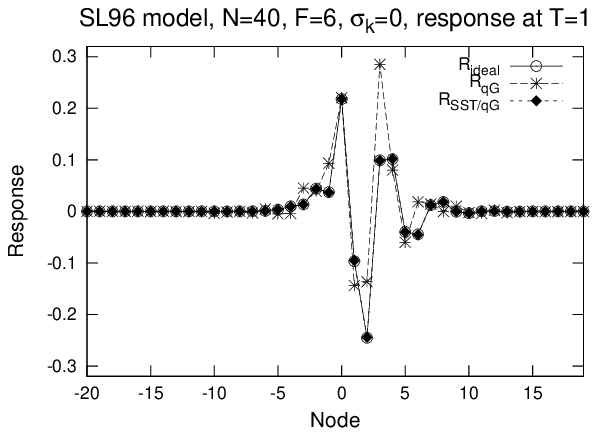}%
\picturehere{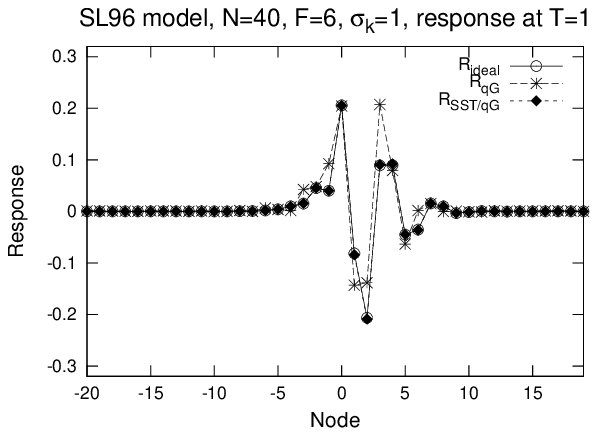}\\%
\picturehere{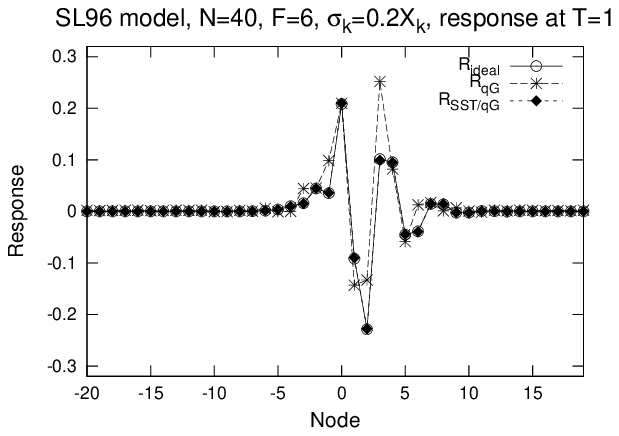}%
\picturehere{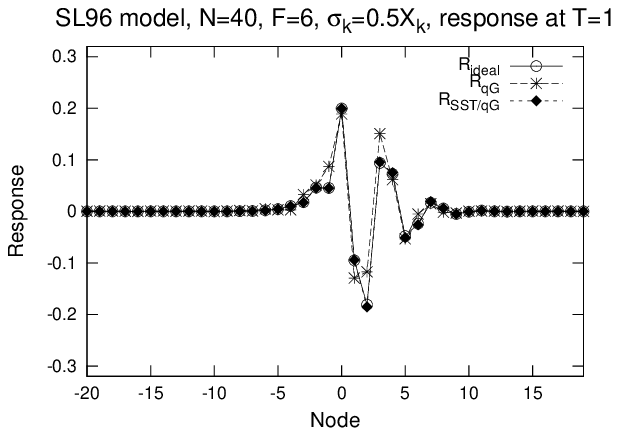}%
\caption{Snapshots of the response operators for SL96 model
at $T=1$, $N=40$, $F=6$.}%
\label{fig:l96_snapshot_t1}%
\end{figure}%
\begin{figure}%
\picturehere{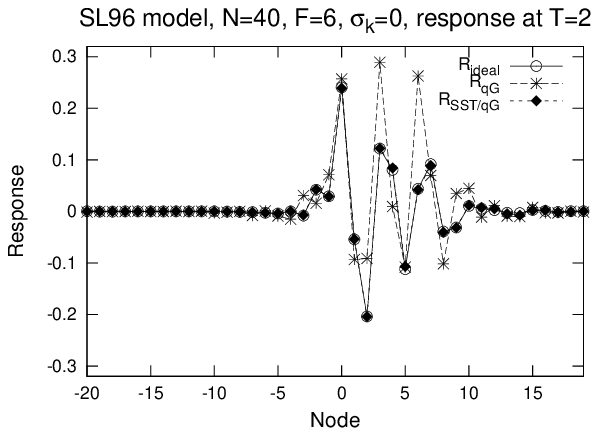}%
\picturehere{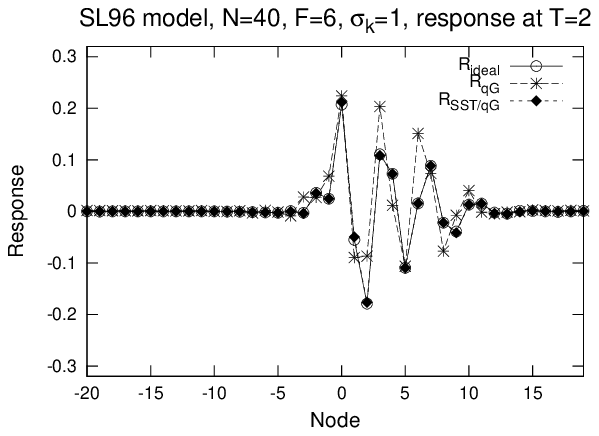}\\%
\picturehere{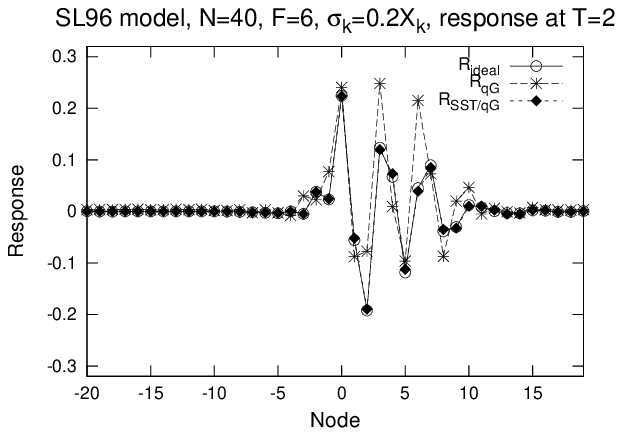}%
\picturehere{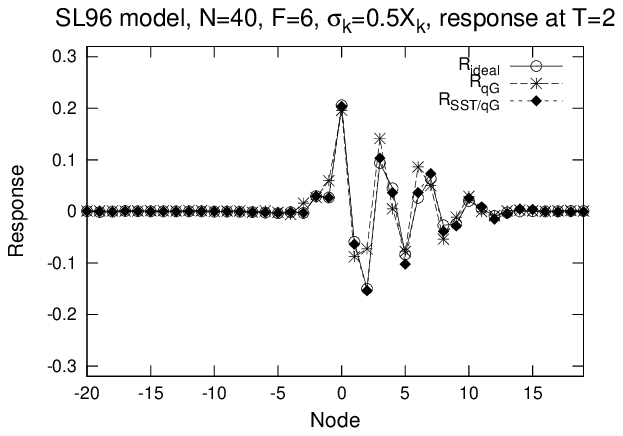}%
\caption{Snapshots of the response operators for SL96 model
at $T=2$, $N=40$, $F=6$.}%
\label{fig:l96_snapshot_t2}%
\end{figure}%
\begin{figure}%
\picturehere{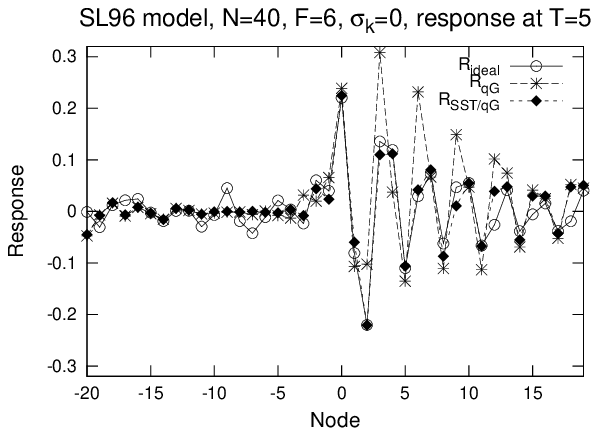}%
\picturehere{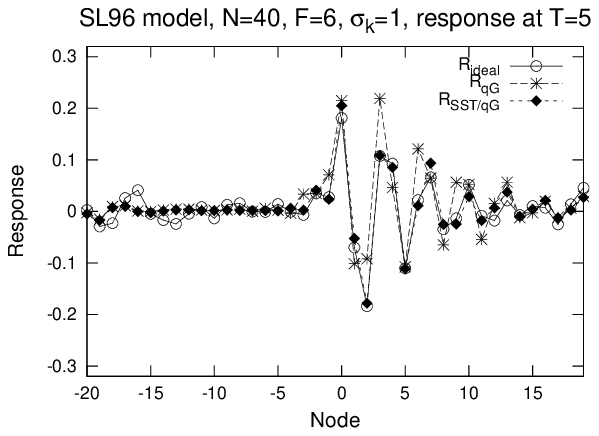}\\%
\picturehere{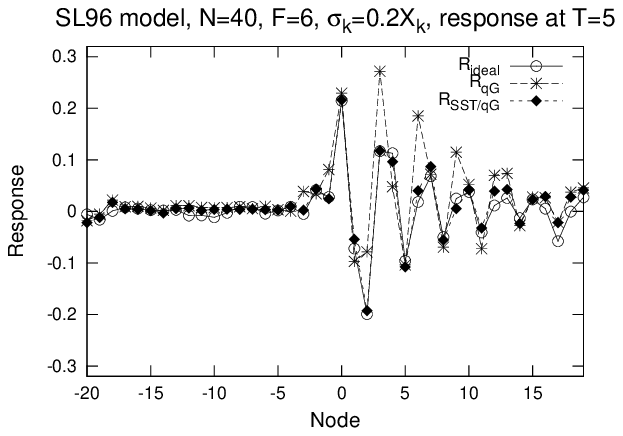}%
\picturehere{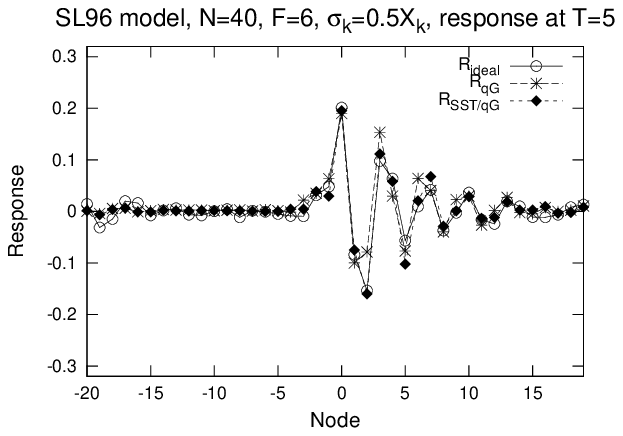}%
\caption{Snapshots of the response operators for SL96 model
at $T=5$, $N=40$, $F=6$.}%
\label{fig:l96_snapshot_t5}%
\end{figure}%
In Figure \ref{fig:l96_errors} we display the $L_2$ relative errors
between the ideal response operator and the FDT response operators,
together with the intrinsic error in the ideal response operator
(which is the result of slight nonlinearity in the ideal response due
to small but finite perturbations). Observe that in the fully
deterministic regime ($F=6$, $\sigma_k=0$) the SST-FDT response
provides a very precise prediction until the time $t\approx 3$, and
then the errors in the SST-FDT grow exponentially rapidly, which is
due to the positive Lyapunov exponents and numerical instability in
the tangent map. On the other hand, the qG-FDT response is not precise
(reaching about 80\% by the time $t=1.5$), due to the fact that the
invariant probability measure associated with the deterministic regime
is highly non-Gaussian, and most probably not continuous with respect
to the Lebesgue measure (that is, it does not even possess a
density). Remarkably, if we look at the stochastically driven regimes
$\sigma_k=1$ (additive noise) and $\sigma_k=0.2X_k$, $\sigma_k=0.5X_k$
(multiplicative noise), we see that the behavior of both the SST-FDT
and qG-FDT responses is qualitatively the same as in the fully
deterministic regime, even though the dynamics is qualitatively
different. Apparently, the level of noise in the two stochastically
driven regimes $\sigma_k=1$ and $\sigma_k=0.2X_k$ is insufficient to
``smooth out'' the invariant probability measure enough for it to
resemble the Gaussian state and to destabilize the computation of the
tangent map. However, in the $\sigma_k=0.5X_k$ multiplicative noise
regime, the errors in the initial qG-FDT response are reduced to about
40\%, which is due to the fact that in this regime the invariant
probability measure is closer to the Gaussian state because of strong
noise. The blended SST/qG-FDT response yields the lowest errors in all
cases, due to its explicit design to avoid numerical instability in
the SST-FDT algorithm.

In Figure \ref{fig:l96_corrs} we show the correlation functions for
the same simulations. Observe that, although significant $L_2$-errors
were observed for the qG-FDT algorithm for the fully deterministic
regime $\sigma_k=0$, its correlations with the ideal response are
generally on the level of around 0.7, which is remarkable. Also, the
correlations of the SST-FDT response with the ideal response are
roughly 1 (nearly perfect correlation) before the numerical
instability manifests itself. As for the blended SST/qG-FDT response,
the best correlations are achieved in the stochastically forced
regimes $\sigma_k=1$ (additive noise) and $\sigma_k=0.2X_k$,
$\sigma_k=0.5X_k$ (multiplicative noise), were the correlations do not
become lower than 0.95 for all response times. For the fully
deterministic case $\sigma_k=0$ the correlations of the blended
SST/qG-FDT response are about 0.8.

In addition to displaying the errors and correlations between the FDT
response operators and the ideal response operator, in Figures
\ref{fig:l96_snapshot_t1}--\ref{fig:l96_snapshot_t5} we show the
instantaneous snapshots of the linear response operators at times
$T=1$, $T=2$ (which are before the SST/qG-FDT cutoff time) and $T=5$
(which is after the SST/qG-FDT cutoff time). Although the linear
response operator at a given time is an $40\times 40$ matrix, it has
the property of translational invariance (just like the L96 model
itself), and, thus, can be averaged along the main diagonal with
wrap-around aliasing of rows (or columns) into a single vector. These
averaged vectors are displayed in Figures
\ref{fig:l96_snapshot_t1}--\ref{fig:l96_snapshot_t5}. Observe that for
the early times of the response $T=1,2$ the SST/qG-FDT response is
virtually indistinguishable from the ideal response. As for the qG-FDT
response, its best performance is observed in the case of strong
multiplicative noise $\sigma_k=0.5X_k$, where the discrepancies
between the qG-FDT and ideal response are not much larger than those
between the SST/qG-FDT response and the ideal response. This is
probably the consequence of the fact that the strong multiplicative
noise changes the invariant probability density of the SL96 model to
the point where it is relatively close to the Gaussian. For other
regimes, by the response time $T=2$ significant errors develop in the
qG-FDT response to the right of the main response diagonal. For the
longer response time $T=5$ and all regimes the blended SST/qG-FDT
response is very similar to the ideal response, while the qG-FDT
response again develops large discrepancies to the right of the main
response diagonal for $\sigma_k=0,1,0.2X_k$. For the strong
multiplicative noise regime, $\sigma_k=0.5X_k$, and response time
$T=5$, the qG-FDT yields lower errors than in the other regimes, but
is still less precise than the SST/qG-FDT response.

\section{Summary}
\label{sec:summary}

The classical fluctuation-dissipation theorem, by its design, is
suitable for computing the linear response for stochastically driven
systems, as it assumes the continuity of the probability measure of
the statistical ensemble distribution with respect to the Lebesgue
measure (which is guaranteed in many stochastically driven
systems). However, the drawback of the classical fluctuation-response
formula is that it requires the probability density together with its
derivative (or their suitable approximations) explicitly in the
response formula. Unfortunately, for complex systems with many
variables such an approximation might not be necessarily available
with required precision.

In this work, we develop the stochastic short-time
fluctuation-dissipation formula (SST-FDT) for stochastically driven
systems which does not require the probability measure of the
statistical state of the system to be known explicitly. This formula
is the analog of the general linear response formula
\cite{AbrMaj4,AbrMaj5,AbrMaj6,EckRue,Rue2} for chaotic (but not
stochastically driven) nonlinear systems. We demonstrate that, before
the numerical instability due to positive Lyapunov exponents occurs,
the SST-FDT for the stochastically driven Lorenz 96 model is generally
superior to the classical FDT formula where the probability density of
the statistical state is approximated by the Gaussian density with the
same mean and covariance (qG-FDT). We test the new SST-FDT formula for
the L96 model with stochastic forcing for both the additive and
multiplicative noise, and observe that the SST-FDT response formula is
generally better than the qG-FDT in both the error and correlation
comparison, before the numerical instability develops in the SST-FDT
response. Additionally, the blended SST/qG-FDT response with a simple
Heaviside blending function clearly performs on top of both the qG-FDT
and SST-FDT in all studied regimes. The results of this work suggest
that the SST/qG-FDT algorithm can be used in practical applications
with stochastic parameterization, such as the climate change
prediction.

\begin{acknowledgment}
The author thanks Ibrahim Fatkullin for helpful comments and
remarks. This work is supported by the NSF CAREER grant DMS-0845760
and the ONR grant N000140610286.
\end{acknowledgment}

\end{document}